\documentclass[11pt]{article}
\usepackage{amssymb}
\usepackage{latexsym}
\usepackage{epsfig}

\usepackage{amstex}
\usepackage{amsxtra}
\usepackage{amstext}
\usepackage{amsfonts}
\usepackage{amsbsy}
\usepackage{amsopn}
\usepackage{amssymb}

\begin{document}
\newtheorem{definition}{Definition}
\newtheorem{theorem}{Theorem}
\newtheorem{fact}{Fact}
\newtheorem{lemma}{Lemma}
\newtheorem{proposition}{Proposition}
\newtheorem{remark}{Remark}
\newtheorem{corollary}{Corollary}
\newtheorem{conjecture}{Conjecture}
\newtheorem{example}{Example}
\input amssym.def
\input amssym
\def\nin{\noindent}

\title{{\bf  Matrix elements of  vertex operators of
              deformed W-algebra and Harish Chandra Solutions to 
Macdonald's  difference equations }}
\author{A. Kazarnovski-Krol } 
\maketitle

\begin{abstract} In this paper we prove that certain matrix elements
of vertex operators of 
deformed W-algebra satisfy Macdonald difference equations and 
form $n!$ -dimensional space of solutions. These solutions are the
analogues of Harish Chandra  solutions with prescribed asymptotic behavior.                 
 We obtain formulas for
analytic continuation as a consequence of braiding properties  of 
 vertex operators of deformed W-algebra.    
\end{abstract}

\section{Introduction}
\label{sec1}
For each $u \in \mathbb C$ and $1 \le i \le n$ 
define   the shift operator  $T_{u,x_i}$ by 
$$T_{u,x_i}f(x_1,x_2,\ldots,x_i,\ldots, x_n)=f(x_1,\ldots, u
x_i, \ldots, x_n) $$
for any complex valued function $f$.
Let  $q, k $ be   real numbers,  $0<q<1$, $0<k<1$. Set also  $t=q^k$. 
Let 
$z=(z_1,z_2,\ldots,z_n)$ be a set of complex variables.
For $m=1,\ldots,n$ consider the following difference operators
\begin{eqnarray}
     D_z^m=D_z^m(q,t)
\nonumber \\= t^{\frac{1}{2}m(m+1)} \sum_{i_1 < i_2 < \ldots < i_m}
\left[ \prod_{ s =1, \ldots, m ; \;  
                  j \notin \{i_1, \ldots, i_m \} }
 \frac{t z_{i_s} - z_j}{z_{i_s}-z_j} \right]  
T_{q, z_{i_1}} T_{q, z_{i_2}} T_{q, z_{i_3}} \ldots T_{q, z_{i_m}}  
\label{f1}.
\end{eqnarray}
Operators $D_z^m$ commute with each other.

For $\gamma=(\gamma_1, \ldots,\gamma_n) \in \mathbb C^n$ let
\begin{equation}
c_{\gamma}^m= \sum_{i_1 < i_2 < \ldots < i_m} \prod_{1 \le s \le m}
q^{\gamma_{i_s}} t^{i_s}
\label{f2}
\end{equation}
Consider the following system of difference equations:
\begin{equation}
D_z^m f(z_1, \ldots, z_n)= c_{\gamma}^m f(z_1, \ldots, z_n),
\quad  \quad m=1,\ldots,n
\label{syst1}
\end{equation}

Define $\delta \in \mathbb C^n$
$$\delta =\left( \frac{n-1}{2}, \frac{n-3}{2}, \ldots, -\frac{n-1}{2}\right).$$
Define $\rho \in \mathbb C^n$ to be
$$\rho=k \delta=
k\left(\frac{n-1}{2}, \frac{n-3}{2}, \ldots, -\frac{n-1}{2}\right)$$

Parameter $\lambda=(\lambda_1, \lambda_2, \ldots, \lambda_n) \in
\mathbb C^n$ is assumed to be {\it GENERIC}
and for the sake of simplicity we impose $\lambda_1+\ldots+\lambda_n=0$.
We assume that for any $l_1,l_2, \ldots, l_n \in \mathbb Z$ ,s.t.
$\sum l_i^2 \ne 0$ the following condition is fullfilled
\begin{equation}
\sum_{i=1}^{n} q^{\lambda_i +l_i} -\sum_{i=1}^n q^{\lambda_i} \ne 0
\label{}
\end{equation}

One can see that
\begin{equation}
c_{\lambda+\rho}^m= \sum_{i_1 < i_2 < \ldots < i_m} \prod_{1 < s < m}
q^{\lambda_{i_s}} t^{\frac{n+1}{2}}.
\label{ }
\end{equation}
Thus for any element $w$ of the symmetric group $S_n$: $w \in S_n$ we have
\begin{equation}
c_{\lambda+\rho}^m= c_{w \lambda+\rho}^m .
\label{f3}
\end{equation}
Here we mean that $w \in S_n$ acts on $\lambda=(\lambda_1, \ldots,
\lambda_n)$ as
$$
w  : \; \; (\lambda_1, \ldots, \lambda_n) \mapsto 
 (\lambda_{w(1)}, \ldots, \lambda_{w(n)})
$$
We consider the following system of difference equations:
\begin{equation}
D_z^m \phi(z) = c_{\lambda+\rho}^m \phi(z)
\label{Mac}
\end{equation}
This  system will be referred to as Macdonald's system of difference
equations.
 
\begin{definition}{(Harish Chandra  solution to Macdonald's system of
difference equations)} Fix $\lambda=(\lambda_1,\ldots, \lambda_n) \in
\mathbb C^n$.
The  solution of Macdonald's system of  difference equations  (\ref{Mac})
in the asymptotic zone $$0<|z_1| < |z_2| < \ldots < |z_n|$$
which has  the following series expansion:
\begin{equation}
\phi(\eta +\rho, z)=z^{\eta+\rho}\sum_{p_1 , p_2 
,\ldots , p_{n-1} \in \Bbb Z_{+}} a(p_1, p_2, \ldots,p_{n-1})
\left(\frac{z_1}{z_2}\right)^{p_1} \ldots
 \left(\frac{z_{n-1}}{z_n}\right)^{p_{n-1}} 
\label{ser2}
\end{equation}
with 
\begin{equation}
a(0,0,\ldots,0)=(-1)^{\frac{(n-1)n}{2}}
\prod_{i <j } q^{\frac{( \eta_i -\eta_j)( \eta_i -\eta_j+k)} {2}}
 \frac{\Gamma_q(1-k)}
{\Gamma_q(\eta_i-\eta_j+1)\Gamma_q(\eta_j-\eta_i +1- k)
}\quad \quad
\label{}
\end{equation}
as the leading asymptotic coefficient,


and with $\eta= w \lambda$, $w \in S_n$,
 will be referred to as Harish Chandra   solution.
The correctness of this definition is not clear apriori and it will be
verified in theorem \ref{th1}. See also corollary \ref{cor1}.
\end{definition}







Let $h_1, h_2, \ldots, h_n$ be the weights of the vector
representation of $sl(n)$. Bosonization of the vertex operators 
$\Phi_{h_i}(z)$ was completed in refs.\cite{AKOS1}, \cite{FF}, \cite{AJMP}, and we recall in 
section \ref{bos}, formula (\ref{vert}) below.
In this paper we prove the following theorem:

{\bf{Theorem.}}
For $\mu=\frac{\lambda}{k}$, $0<k<1$ the following
matrix elements 
\begin{eqnarray} \prod z_i^{   k (\omega_1,\delta)
-\frac{k}{2}(\omega_1, \omega_1)}
 \prod_{i < j} g_1\left(\frac{z_i}{z_j}\right)^{-1}
\prod z_i^{-(i-1)k(\omega_1, \omega_1)} \nonumber \\
\times <0,\mu^*| \Phi_{h_{w(n)}}(z_n) \ldots \Phi_{h_{w(2)}}(z_2)  
\Phi_{h_{w(1)}}(z_1)|0,\mu>
\end{eqnarray}
are the Harish Chandra  solutions $\phi(w \lambda+\rho,z)$ 
to Macdonald's difference equations. The Harish Chandra solution $\phi(w
\lambda+\rho,z)$ converges to analytic function for sure at least  for $$0<|z_1| < |z_2|
< \ldots < |z_n|.$$
see theorem \ref{th1} below.
\begin{remark}
Harish Chandra solution $\phi(w \lambda+\rho, z)$ is also an
eigenfunction of difference operator
$D^1_z( \frac{1}{q}, \frac{1}{t})$ as it immediately follows from
\begin{eqnarray}
 D^{n-1}_z(q,t) \phi(w \lambda+\rho,z) \nonumber \\ =
t^{\frac{(n-1)n}{2} } \sum_{i_1 <i_2 < \ldots < i_{n-1}} \left[
\prod_{s=1, \ldots, n-1, \; j \notin \{i_1, \ldots, i_{n-1}\}
} \frac{t z_{i_s}-z_j}{z_{i_s}-z_j}  \right] 
T_{q,z_{i_1}} \ldots T_{q, z_{i_{n-1}}
} \phi(w \lambda +\rho, z)\nonumber \\=
t^{\frac{n(n-1)}{2} } t^{n-1} \sum_{j}
\prod_{i_s \ne j  } \frac{ \frac{1}{t}z_j - z_{i_s}} 
{z_j-z_{i_s}} T_{q^{-1},z_j}
\phi(w \lambda+\rho, z) \nonumber \\ =
t^{\frac{n(n+1)}{2}} D^1_z(\frac{1}{t}, \frac{1}{q})
\phi( w \lambda + \rho, z)
\end{eqnarray}

\end{remark}

\begin{remark}
The system of difference equations (\ref{syst1}) is a generalization of 
radial parts of Laplace-Casimir operators on Riemannian symmetric
space $G/K$ and interpolates between real and p-adic case.            
\end{remark}

Before proceeding further we would like to consider the case of root
system of type $A_1$ in more detail.

\section{$A_1$-case}
\label{sec2}
\subsection{}
In this section we consider the case of root system $A_1$ as a basic
example. It is also used as a base of induction in the proof of
theorem \ref{th1} below.

There is only one difference equation:
\begin{equation}
t \left(\frac{t z_1 -z_2}{z_1-z_2} T_{q,z_1} + \frac{t z_2-z_1}{z_2-z_1}
T_{q,z_2} \right) f(z_1,z_2)=
c^1_{\lambda+\rho} f(z_1,z_2)
\label{f4}
\end{equation}
We are looking for the solution of the equation  (\ref{f4}) of the form:
\begin{equation}
f(z_1,z_2)=z_1^{\lambda_1+\frac{k}{2}} z_2^{\lambda_2-\frac{k}{2}} \sum_{n=0}^{\infty} a(n)
\left(\frac{z_1}{z_2}\right)^{n},
\label{f5} 
\end{equation}
with $a(0)=1$.

$$t=q^k$$

Then 
$$
c_{\lambda+\rho}^1=t(q^{\lambda_1+\frac{k}{2}} + t
q^{\lambda_2-\frac{k}{2}})=
t^{\frac{3}{2}} (q^{\lambda_1}+q^{\lambda_2})
$$




The solution is given by q-hypergeometric function:
\begin{equation}
f(z_1,z_2)=z_1^{\lambda_1+\frac{k}{2}} z_2^{\lambda_2-\frac{k}{2}}
 {F_q}(k, \lambda_1-\lambda_2+k, \lambda_1-\lambda_2 +1, \frac{q}{t}
\frac{z_1}{ z_2})
\label{f6}
\end{equation}

$$
F_q(a,b,c,z)=1+\sum_{n=1}^{\infty} \frac{(1-q^a)\ldots (1-q^{a+n-1})
(1-q^b) \ldots (1-q^{b+n-1})}
{(1-q) \ldots(1-q^{n}) (1-q^c) \ldots (1-q^{c+n-1})} z^n
$$

Note that if we let $\lambda_1-\lambda_2+k=-n$ to be negative
integral, then hypergeometric function terminates and 
$$
z_2^n F_q( k, -n, -n-k+1,  \frac{q}{t}
\frac{z_1}{z_2})
$$
provides Macdonald polynomial of $2$ variables of degree $n$.

Examples.

0. $n=0$
$$
z_2^0 F_q(k, 0, 1-k, \frac{q}{t}\frac{z_1}{z_2})= 1
$$

1. $n=1$
$$
z_2 F_q(k, -1, -k,   \frac{q}{t}
\frac{z_1}{z_2})=
z_2\left(1 +\frac{q}{t} \frac{(1-q^k)(1-q^{-1})}{(1-q)(1-q^{-k})}
\frac{z_1}{z_2}\right)=
z_2+z_1
$$

2. $n=2$
\begin{eqnarray}
z_2^2 F_q(k, -2, -1-k, \frac{q}{t} \frac{z_1}{z_2})=
z_2^2\left(1 + \frac{(1-q^k)(1-q^{-2})}{(1-q)(1-q^{-1-k})} \frac{q}{t}
\frac{z_1}{z_2}+\right.\nonumber \\ \left.
\frac{(1-q^k)(1-q^{k+1})(1-q^{-2})(1-q^{-1})}{(1-q)(1-q^2)(1-q^{-1-k})(1-q^{-k})}
(\frac{q}{t})^2 (\frac{z_1}{z_2})^2\right)\nonumber \\=
z_2^2 + z_1^2 + \frac{(1-q^k)(1-q^{-2})}{(1-q)(1-q^{-1-k})}
\frac{q}{t} z_1 z_2\nonumber \\=
z_2^2 + z_1^2+ \frac{(1-t)(1+q)}{(1-tq)} z_1 z_2
\label{}
\end{eqnarray}

3. $n=3$.
\begin{eqnarray}
z_2^3 F_q(k, -3, -2-k, \frac{q}{t} \frac{z_1}{z_2})=
z_2^3  \left( 1 + \frac{(1-q^k)( 1-q^{-3})}{(1-q)(1-q^{-2-k})} \frac{q}{t}
\frac{z_1}{z_2}\right.\nonumber \\
+\frac{(1-q^k)(1- q^{k+1})( 1-q^{-3})(1-q^{-2})}
{(1-q)(1-q^2)(1-q^{-2-k})(1-q^{-1-k})}
(\frac{q}{t} \frac{z_1}{z_2})^2\nonumber \\
\left.+
\frac{(1-q^k)(1- q^{k+1})(1-q^{k+2})( 1-q^{-3})(1-q^{-2})(1-q^{-1})}
{(1-q)(1-q^2)(1-q^3)(1-q^{-2-k})(-q^{-1-k})(1-q^{-k})}
(\frac{q}{t} \frac{z_1}{z_2})^3  \right) \nonumber \\
=
z_2^3 +z_1^3 + (z_1^2 z_2 +z_2^2 z_1) \frac{(1-t)(1+q+q^2)}{(1-q^2 t)}
\end{eqnarray}

4. $n=4$
\begin{eqnarray}
z_2^4 F_q(k, -4, -3-k, \frac{q}{t} \frac{z_1}{z_2})=
(z_2^4 +z_1^4) + \frac{(1-t)(1+q+q^2+q^3)}{(1-q^3 t)}(z_2^3 z_1 +z_1^3 z_2)
\nonumber \\+ \frac{(1+q^2)(1+q+q^2)(1-t)(1-qt)}{(1-q^2 t)(1-q^3 t)}z_2^2 z_1^2
\end{eqnarray}

One can see in these examples that one gets  indeed
symmetric polynomials.

\smallskip

\subsection{}

We will use the standard notations for  $q-\Gamma$ -function and
$\Theta$ function:
$$
(z;q)_{\infty}=\prod_{i=0}^\infty(1-z q^i)
$$

$$
\Gamma_q(x)=\frac{(q;q)_{\infty} (1-q)^{1-x}}
{(q^x;q)_{\infty}}
$$

$$
\Theta_q(z)=(z;q)_{\infty}(q z^{-1};q)_{\infty}(q;q)_{\infty}
$$
The subscript $q$ in notations of $\Theta$-function will be omitted
sometimes.




There is the following integral representation  for q-hypergeometric
function with usual integration:

\begin{eqnarray}
\int \frac{dy}{2 \pi i y} y^{\lambda_1-\lambda_2} 
\frac{ (q^{\frac{1+k}{2}} y^{-1};q)_{\infty}
(q^{\frac{1+k}{2}} y z  ;q)_{\infty}}
{(q^{\frac{1-k}{2}} y^{-1};q)_{\infty} (q^{\frac{1-k}{2}} y z;
q)_{\infty}} \nonumber \\
=
q^{\frac{(\lambda_1-\lambda_2)(1-k)}{2}} \frac{\Gamma_q(\lambda_1-\lambda_2+k)}{\Gamma_q(\lambda_1-\lambda_2+1)\Gamma_q(k)}
F_q(\lambda_1-\lambda_2+k, k, \lambda_1-\lambda_2+1, q^{1-k} z)
\label{f7}
\end{eqnarray}
Here the contour of  integration for $y$ starts and ends at $0$ and
encloses poles $y=q^{\frac{1-k}{2}+n}$  .
The integrand of (\ref{f7}) has two series of poles:
$y=q^{\frac{1-k}{2}+n}$ nd $ y= q^{\frac{k-1}{2}-n}z^{-1}$ . They do not
overlap if
$|q^{\frac{1-k}{2}}|< |q^{\frac{k-1}{2}} z^{-1}|$ or
$|z|<q^{k-1}$. Since we assume that $ 0 < k < 1$ it for sure converges
for $ |z| <1$.
The importance of this integral representation for q-hypergeometric
function with usual integration was emphasized in ref. \cite{LP2}.

For the purposes of completeness we recall the famous Ramanujan's integral 

\begin{equation}
\int_0^{\infty} x^{c-1} \frac{(-ax; q)_{\infty}}{(-x;q)_{\infty}} dx =
\frac{\pi}{\sin \pi c} \frac{(q^{1-c};q)_{\infty} (a;q)_{\infty}}
{(q;q)_{\infty}(a q^{-c};q)_{\infty}}
\label{}
\end{equation}

Also one has:

\begin{eqnarray}
(z_1 z_2)^{\lambda_1+\frac{k}{2}}
\int \frac{dy}{2 \pi i y} y^{\lambda_2-\lambda_1}
\frac{(q^{\frac{1+k}{2} } \frac{z_1}{y};q)_{\infty}
(q^{\frac{1+k}{2}} \frac {y}{z_2};q)_{\infty}}
{(q^{\frac{1-k}{2}} \frac{z_1}{y};q)_{\infty}
(q^{\frac{1-k}{2}} \frac{y}{z_2} ;q)_{\infty}}\nonumber \\
=
q^{\frac{(\lambda_2-\lambda_1)(1-k)}{2}}
\frac{\Gamma_q(\lambda_2-\lambda_1+k)}
{\Gamma_q(k) \Gamma_q(\lambda_2-\lambda_1+1)} 
z_1^{\lambda_2+\frac{k}{2}} z_2^{\lambda_1-\frac{k}{2}}
F_q(k, \lambda_2-\lambda_1+k, \lambda_2-\lambda_1+1, q^{1-k} \frac{z_1}{z_2})
\label{f8}
\end{eqnarray}



Note that $y^{\lambda_2-\lambda_1}$ can be replaced by appropriate
ratio of $\Theta$-functions:

\begin{eqnarray}
\int \frac{dy}{2 \pi i y}
\frac{\Theta_q(q^{\lambda_2-\lambda_1 +\frac{1+k}{2}} \frac{z_1}{y})}
{\Theta_q(q^{\frac{1+k}{2}} \frac{z_1}{y})}
\frac{(q^{\frac{1+k}{2}} \frac{z_1}{y};q)_{\infty}}
{(q^{\frac{1-k}{2}} \frac{z_1}{y};q)_{\infty}}
\frac{(q^{\frac{1+k}{2}} \frac{y}{z_2};q)_{\infty}}
{(q^{\frac{1-k}{2}} \frac{y}{z_2};q)_{\infty}}=
\nonumber \\
\frac{\Theta_q(q^{\lambda_2-\lambda_1+k})}
{\Theta_q(q^k)} 
\frac{\Gamma_q(\lambda_2-\lambda_1+k)}
{\Gamma_q(k) \Gamma_q(\lambda_2-\lambda_1+1)}
F_q(k, \lambda_2-\lambda_1+k, \lambda_2-\lambda_1+1, q^{1-k}
\frac{z_1}{z_2}) \nonumber \\
=\frac{\Gamma_q(1-k)}{\Gamma_q( \lambda_1-\lambda_2+1-k)
\Gamma_q(\lambda_2-\lambda_1+1)}
F_q(k, \lambda_2-\lambda_1+k, \lambda_2-\lambda_1+1, q^{1-k} \frac{z_1}{z_2})
\label{ff}
\end{eqnarray}

Here the integrand is single-valued and the contour of integration can
be  taken away from the origin.

\begin{eqnarray}
(z_1 z_2)^
{\lambda_1+\frac{k}{2}} 
\int y^{\lambda_2-\lambda_1+k}z_1^{-k} 
\frac{(q^{\frac{1+k}{2}}
\frac{y}{z_1};q)_{\infty}}
{(q^{\frac{1-k}{2}}
\frac{y}{z_1};q)_{\infty}}
\frac{(q^{\frac{1+k}{2}}\frac{y}{z_2};q)_{\infty}}
{(q^{\frac{1-k}{2}}\frac{y}{z_2};q)_{\infty}}
z_2^{-k}
\frac{\Theta_q(q^{\lambda_1-\lambda_2+\frac{1-k}{2}}
\frac{y}{z_1} )}
{\Theta_q(q^{\frac{1+k}{2}}
\frac{y}{z_1})} 
\left(\frac{z_1}{y}\right)^{\lambda_2-\lambda_1+k} \frac{dy}{2 \pi i y}\nonumber \\=
\frac{\Gamma_q(1-k)}{\Gamma_q(\lambda_1-\lambda_2+1-k) \Gamma_q(\lambda_2-\lambda_1+1)} 
z_1^{\lambda_2+\frac{k}{2}}
z_2^{\lambda_1-\frac{k}{2}}
F_q(k, \lambda_2-\lambda_1+k, \lambda_2-\lambda_1+1, q^{1-k} \frac{z_1}{z_2}
)
\end{eqnarray}

Here the contour of integration is chosen to be closed curve around
the origin  
$$q^{\frac{1-k}{2}}|z_1| < |y| < |z_1|$$ 
enclosing the poles 
$$y=q^{\frac{1-k}{2}+n} z_1$$
$n=0,1,2,\ldots$

\smallskip
Also, one has
\begin{eqnarray}
(z_1 z_2)^{\lambda_1+\frac{k}{2}}
\int y^{\lambda_2-\lambda_1+k} y^{-k} 
\frac{(q^{\frac{1+k}{2}} \frac{z_1}{y};q)_{\infty}}
{(q^{\frac{1-k}{2}} \frac{z_1}{y};q)_{\infty}}
z_2^{-k} \frac{(q^{\frac{1+k}{2}} \frac{y}{z_2};q)_{\infty}}
{(q^{\frac{1-k}{2}} \frac{y}{z_2};q)_{\infty}} 
\frac{\Theta_q(q^{\lambda_1-\lambda_2+\frac{1+k}{2}} \frac{y}{z_2})}
{\Theta_q(q^{\frac{1+k}{2}} \frac{y}{z_2})}
\left(\frac{z_2}{y}   \right)^{\lambda_2-\lambda_1} \frac{d y}{ 2 \pi
i y}\nonumber \\
=
\frac{\Gamma_q(1-k)}{\Gamma_q(\lambda_2-\lambda_1+1-k)\Gamma_q(\lambda_1-\lambda_2
+1)}
  z_1^{\lambda_1+\frac{k}{2}} z_2^{\lambda_2-\frac{k}{2}}
F_q(k, \lambda_1-\lambda_2+k, \lambda_1-\lambda_2+1, q^{1-k}
\frac{z_1}{z_2})
\end{eqnarray}

Here the contour of integration is chosen to be closed curve around
the origin 
$$ q^{\frac{1-k}{2}} |z_2| < |y| < |z_2|$$
enclosing the two series of poles:
$$
y=q^{\frac{1-k}{2}+n} z_1$$
$$
y=q^{\frac{1-k}{2}+n} z_2$$
$n=0,1,2,\ldots$
 





\subsection{}
For the q-hypergeometric function analytic continuation is given:
\begin{eqnarray}
F_q(a,b,c; z)=
\frac{\Gamma_q(c) \Gamma_q(b-a)}{\Gamma_q(b) \Gamma_q(c-a)}
\frac{\Theta_q(q^a z)}{\Theta_q(z)}
F_q(a, a-c+1, a-b+1, q^{c+1-a-b} z^{-1}) \nonumber \\+
\frac{\Gamma_q(c) \Gamma_q(a-b)}{\Gamma_q(a) \Gamma_q(c-b)}
\frac{\Theta_q(q^b z)}{\Theta_q(z)}
F_q(b, b-c+1, b-a+1, q^{c+1-a-b} z^{-1})
\label{f9}
\end{eqnarray}




The series for $F_q(k, \lambda_1-\lambda_2+k, \lambda_1-\lambda_2+1,
q^{1-k}z)$
converges if $|q^{1-k}z| <1$,  equivalently, $|z|< q^{k-1}$.






If we normalize 
the two solutions to difference   
equation (\ref{f4}) as:

\begin{eqnarray}
\hat f ( k, \lambda_1-\lambda_2+k, \lambda_1-\lambda_2+1, 
{z_1},{z_2})\nonumber \\=
\frac{\Gamma_q(1-k)}{\Gamma_q(\lambda_2-\lambda_1+1-k)
\Gamma_q(\lambda_1-\lambda_2+1)}
z_1^{\lambda_1+\frac{k}{2}} z_2^{\lambda_2-\frac{k}{2}}
 F_q(k, \lambda_1-\lambda_2+k,
\lambda_1-\lambda_2+1, q^{1-k} \frac{z_1}{z_2})
\end{eqnarray}

and analogously
\begin{eqnarray}
\hat f(k,
 \lambda_2-\lambda_1+k, \lambda_2-\lambda_1+1,  
 {z_1},{z_2})\nonumber \\=
\frac{\Gamma_q(1-k)}{ \Gamma_q(\lambda_2-\lambda_1+1)
\Gamma_q(\lambda_1-\lambda_2+1-k)} 
z_1^{\lambda_2+\frac{k}{2}} z_2^{\lambda_1-\frac{k}{2}}
F_q( k, \lambda_2-\lambda_1+k,
\lambda_2-\lambda_1+1, q^{1-k} \frac{z_1}{z_2})
\end{eqnarray}

Then formulas for analytic continuation are written as:
\begin{eqnarray}
\hat f(k, \lambda_1-\lambda_2+k,\lambda_1-\lambda_2 +1, 
z_1,z_2)\nonumber \\=
\left(\frac{z_1}{z_2}
\right)^k
\frac{ \Theta_q(q^{\lambda_1-\lambda_2+k}) \Theta_q(
\frac{z_2}{z_1})}
{\Theta_q(q^{\lambda_1-\lambda_2}) \Theta_q(q^k \frac{z_2}{z_1})}
\hat f(k, \lambda_2-\lambda_1+1, z_2, z_1)\nonumber \\+
\left(\frac{z_1}{z_2}\right)^{\lambda_1-\lambda_2+k}
\frac{ \Theta_q(q^k) \Theta_q( q^{\lambda_2-\lambda_1} \frac{z_2}{z_1})}
{\Theta_q(q^{\lambda_2-\lambda_1}) 
\Theta_q(q^k  \frac{z_2} {z_1})}
\hat f( k, \lambda_1-\lambda_2+k, \lambda_1-\lambda_2+1, 
z_2, z_1)
\end{eqnarray}

If we assume, that $ 0 <k <\frac{ 1}{2}$ , then     for $z_1=1$,
$z_2=t=q^k$ ,$\frac{z_1}{z_2}=q^{-k}$, 
we have 
$$
F_q(k, \lambda_1-\lambda_2+k, \lambda_1-\lambda_2+1, q^{1-2k})=
\frac{ \Gamma_q(\lambda_1-\lambda_2+1) \Gamma_q(1-2 k)}
{\Gamma_q(\lambda_1-\lambda_2+1-k) \Gamma_q(1-k)}
$$



\section{ Vertex operators of deformed W-algebra}
\label{sec3}
The deformed $W$-algebras is a new, interesting subject. We refer the
reader to refs.\cite{AKOS1},\cite{FF}  , \cite{AJMP}, 
\cite{AKOS2},  \cite{JLMP},\cite{LL}, \cite{LP}, \cite{LP2}.
In this section we closely follow to the works  \cite{AJMP},\cite{LP}.

\subsection{ Bosons.}
\label{bos} 
Let $x$ be a real parameter, $0<x<1$ , $r>1$.
Consider the bosonic oscillators $\beta_m^j$ 
$$(1 \le j \le n-1, m \in \mathbb Z \setminus \{0\})$$
with commutation relations 
\begin{equation}
[{\beta_m^j} , {{\beta_{m^{\prime}}}^p }]=
\left\{ \begin{array}{lll}
m \frac{x^{(n-1)m}- x^{-(n-1)m} }{x^{nm}-x^{-nm}}
 \frac{x^{(r-1)m} -x^{-(r-1)m}}{x^{rm}-x^{-rm}}
\delta_{m+m^{\prime},0} &\quad 
if \quad  (j=p) \quad  \\ 
-m x^{sgn(j-p)nm } \frac{x^m-x^{-m}}
 {x^{nm}-x^{-nm}} \frac{x^{(r-1)m}-x^{-(r-1)m}}{x^{rm}-x^{-rm}} 
\delta_{m+m^{\prime},0}&
\quad  if \quad (j \ne p)\end{array}
\right. 
\label{f12}
\end{equation}

Define
$\beta_m^n$ by 
$$\sum_{j=1}^n x^{-2 jm} \beta_m^j=0.$$
These oscillators were introduced in refs. \cite{AKOS1}
, \cite{FF}.

Let $\alpha_1, \alpha_2, \ldots, \alpha_{n-1}$ be the simple roots of
root system of type $A_{n-1}$, $\omega_1, \ldots, \omega_{n-1}$ are
fundamental weights, i.e. $<\alpha_j,\omega_i>=\delta_{ij}$,
$\Sigma_{+}$ denotes the  set of positive roots.

The zero mode operators $P_{\alpha}, Q_{\alpha}$ indexed by $\alpha
\in P = \oplus\mathbb Z \omega_i$ are by definition $\mathbb Z$ linear
in $\alpha$ and satisfy
\begin{equation}
[i P_{\alpha}, Q_{\beta}]=<\alpha, \beta> \qquad\qquad  (\alpha, \beta
\in P)
\label{eQ1}
\end{equation}

Define the Fock module $\mathcal F_{l,s}$ generated by $\beta_{-m}^j$ , $m>0$
with  the highest weight vector $|l,s>$ :
\begin{eqnarray}
 P_{\alpha} |l,s>= <\alpha, \sqrt{\frac{r}{r-1}}l - \sqrt{\frac{r-1}{r}}s
>|l,s>,
\nonumber \\
\beta_m^j |l,s>=0, \quad {\rm for} \quad m>0 \nonumber \\
|l,s>= e^{i \sqrt{\frac{r}{r-1}}Q_l- i \sqrt{\frac{r-1}{r} }Q_s} |0,0>
\end{eqnarray}

Expand also (\ref{eQ1}) linearly on $\lambda$ ($\mu$)  
\begin{equation}
[i P_{\alpha}, Q_{\lambda}]=<\alpha, \lambda> \quad  {\rm{for}} \quad 
 (\alpha \in P)
\end{equation}

And in general define state $|\lambda, \mu>$ as 
\begin{eqnarray}
|\lambda, \mu> = e^{i {\sqrt{\frac{r}{r-1}}}Q_{\lambda} -
i {\sqrt{\frac{r-1}{r}}Q_{\mu}}} 
| 0,0 > 
\label{}
\end{eqnarray}

As it is proven in ref. \cite{FF} the notion of Verma module and PBW basis
 make sense for the  deformed (quantum)W-algebra.
We define $< (\lambda,\mu)^*|$ as the  functional, such that  
 $< (\lambda,\mu)^*|| \lambda, \mu>=1$
and $0$ on other subspaces according to PBW basis.

The basic operators $U_{-\alpha_j}(z), U_{\omega_j}(z)$ are defined as follows:
\begin{eqnarray}
U_{-\alpha_j}(z) = e^{i \sqrt{\frac{r-1}{r}}(Q_{\alpha_j}- i
P_{\alpha_j} \log z)} :e^{\sum_{m \ne 0} \frac{1}{m}(\beta_m^j -
\beta_m^{j+1})(x^j z)^{-m}}: \\
U_{\omega_j}(z)=e^{-i \sqrt{\frac{r-1}{r}}(Q_{\omega_j} -i P_{\omega_j}
\log z)}:
e^{-\sum_{m \ne 0} \frac{1}{m}\sum_{p=1}^j x^{(j-2 p+1)m} \beta_m^p
z^{-m}}: 
\label{sd1} 
\end{eqnarray}
Note that
\begin{equation}
e^{i \sqrt{\frac{
r-1}{r}}(Q_{\beta}-
i P_{\beta} \log z)}=
z^{\frac{r-1}{2r}(\beta,\beta)}
e^{i\sqrt{\frac{r-1}{r}}
Q_{\beta}} e^{\sqrt{\frac{r-1}{
r}} P_{\beta} \log z}
\label{feq1}
\end{equation}
as it follows from Campbell-Hausdorf formula. 
The normal ordereing $::$ is defined such that $\beta^j_m $ stands to
the right of $\beta^l_{-m}$ for $m>0$ and symbols $P$ stand to the
right of symbols $Q$.




\subsection{Vertex operators.}
We set
$$
[v]=x^{\frac{v^2}{r}-v} \Theta_{x^{2r}}(x^{2v})
$$
and
\begin{eqnarray}
f(v,w)= \frac{[v+\frac{1}{2} -w]} {[v-\frac{1}{2}]} 
\label{f14}
\end{eqnarray}
For $j<i$ let $\pi_{ji}$ act on $\mathcal F_{l,s}$ as 
$<\alpha_j+\ldots+\alpha_{i-1}  , r l - (r-1) s>$, i.e.
\begin{equation}
\pi_{ji}\Big|_{\mathcal F_{l,s}}=<\alpha_j+\ldots+\alpha_{i-1}  , r l - (r-1) s>
\label{al1}
\end{equation}

Let $h_1, \ldots , h_n$ be  the weights of the vector
representation of $sl(n)$,i.e. $h_1=\omega_1$, $h_2=h_1-\alpha_1$, 
$h_3=h_1-\alpha_1-\alpha_2$, ...,
$h_n=h_1-\alpha_1-\ldots -\alpha_{n-1}$. Note that $h_1+h_2+\ldots+h_n=0$.

 Vertex operators  $\Phi_{h_i}(z)$ of the first type 
are defined as follows:
\begin{eqnarray}
\label{vert}
\Phi_{h_i}(z)=
\oint \prod_{j=1}^{i-1} \frac{d y_j}{2 \pi  i y_j}
U_{\omega_1}(z) U_{-\alpha_1}(y_1) U_{-\alpha_2}(y_2)
\ldots U_{-\alpha_{i-1}}(y_{i-1})
\prod_{j=1}^{i-1} f(v_j-v_{j-1}, \pi_{ji})
\label{f15}
\end{eqnarray}
cf. \cite{AJMP}, see also \cite{LP}, \cite{JLMP}.
Here $x^{2 v_l}=y_l$, $l=1, \ldots, i-1$,  $x^{2 v_0}=z$.

The phase factor  $\prod f(...)$ simultaneously 
makes the integrand single-valued and produces
the  leading asymptotic exponent of $z$.
The contours for integration are chosen to be simple closed curves
around the origin,
s.t. $$ x |z |< |y_1| < | z | $$
$$x |y_1| < |y_2| < |y_1| $$
and so on.

\begin{figure}
\epsfig{file=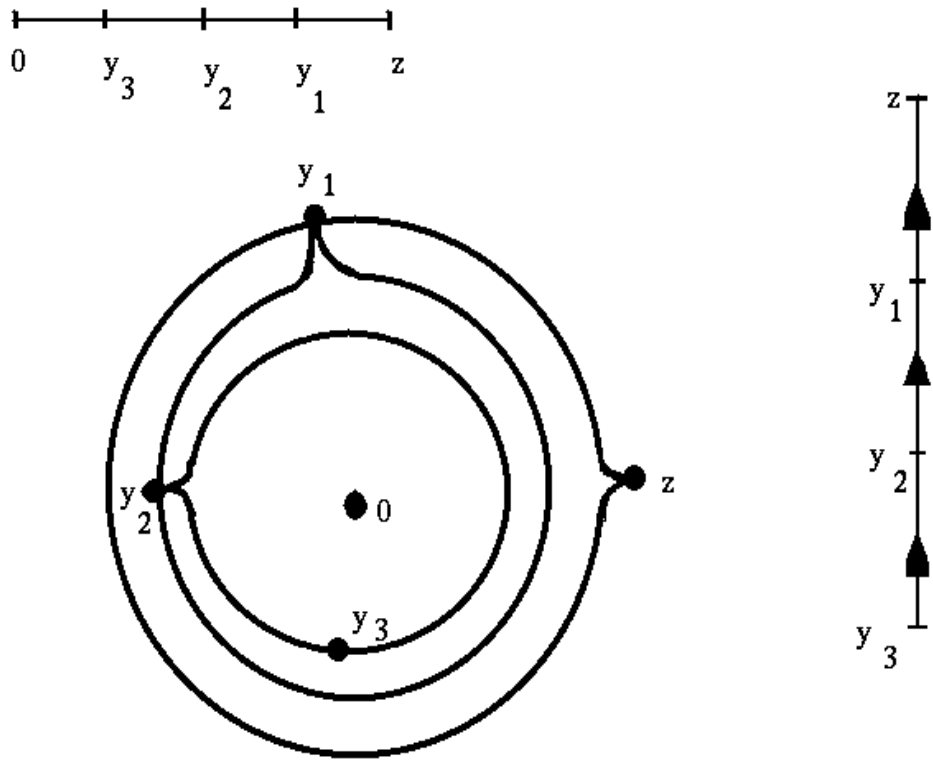, angle=0, height=11cm, width=14cm }
\caption{Cycle for integration for  vertex operator
associated with vector representation of $sl(n)$
(nondeformed case). Note that all internal contours of integration 
are movable.}
\label{125f2}
\end{figure}

In the limit $q \mapsto 1$ the contours of integration tend to those
indicated on fig. \ref{125f2}. Thus the phase factor  $\prod f (...)$
might be considered an analogue of these contours of integration. The
phase factor was introduced in \cite{LP} in the case of deformed
Virasoro, and developed in \cite{JLMP}, \cite{AJMP} and is a very
important ingredient of the construction.
The vertex operators are the same as \cite{AKOS1}, \cite{FF} up to the phase
factor $\prod f(...)$ and up to the shift of parameter of Fock space
by $-\delta$. This shift is caused by  (\ref{sd1}), (\ref{feq1})     and is
designed  in order to alleviate formulas (\ref{f14}),(\ref{al1}),
(\ref{f15}), as well as          (\ref{for1}), (\ref{f16}) ,
and (\ref{f17}) below. It is similar and related  to appearance of $\rho$ in
formulas in harmonic analysis.

Set 
\begin{equation}
\mu_{ij}=(\mu_i-\mu_j)(r-1)
\label{for1}
\end{equation}


Boltzmann weights:
\begin{equation}
W \left[\left.   \begin{array}{lll}
\mu + h_{i}+h_{j} \quad \mu+h_i \\
\mu+h_j \quad\qquad \qquad  \mu \\
\end{array} \right| z
\right]=
r_1(v) 
\frac{ [v][\mu_{ij}-1]}{[v-1] [\mu_{ij}]}
%
%
\label{f16}
\end{equation}

\begin{equation}
W \left[\left. \begin{array}{lll}
\mu+h_i+h_j \quad  \mu+h_j \\
\mu+h_j\quad \quad \qquad  \mu \\
\end{array} \right| z \right]
=
r_1(v) \frac{[v-\mu_{ij}][1]}
{[v-1] [\mu_{ij}]}
%
\label{f17}
\end{equation}

\begin{equation}
r_1(v)=z^{\frac{r-1}{r}\frac{n-1}{n}}\frac{g_1(z^{-1})}{g_1(z)} \quad
(z=x^{2v})
\label{f18}
\end{equation}
\begin{equation}
g_m(z)=\frac{ \left\{x^{m+1} z \right\} \left\{x^{2r+2n-m-1} z \right\}}
{ \left\{x^{2r +m-1} z \right\} \left\{x^{2n-m+1} z \right\}}
\label{g1} 
\end{equation}

$$
\left\{z\right\}=(z;x^{2r},x^{2n})_{\infty}=
\prod_{i_1,i_2=0}^{\infty}( 1-(x^{2r })^{i_1} (x^{2n})^{i_2} z)
$$

For $i \ne j$  vertex operators $\Phi_{h_i}(z_1)$,
$\Phi_{h_j}(z_2)$
 satisfy the following commutation relations:
\begin{eqnarray}
\Phi_{h_i}(z_1) \Phi_{h_j}(z_2)=
W\left[\left. \begin{array}{lll}\mu +h_j +h_i  \quad \mu +h_j \\
\mu +h_j  \qquad \qquad \mu \end{array}\right| \frac{z_2}{z_1} \right]  \Phi_{h_i}(z_2)
\Phi_{h_j}(z_1) \nonumber \\+
W\left[ \left. \begin{array}{lll}
\mu +h_j +h_i \quad \mu +h_i \\
\mu +h_j \qquad \qquad \mu \end{array}\right| \frac{z_2}{z_1}  \right]  
\Phi_{h_j}(z_2) \Phi_{h_i}(z_1)
\label{f19}
\end{eqnarray}

Here it is assumed that  braiding is realized from the asymptotic zone 
$|z_1| > |z_2|$ to the asymptotic zone $|z_2| > |z_1|$.

Bosonization of relation (\ref{f19}) was initiated in ref. \cite{LP2},
and developed in \cite{LP}, \cite{AJMP}.
The relation (\ref{f19}) is in agreement with proposition \ref{prop4}
below.

\section{Matrix elements of   vertex operators of
 deformed W-algebra provide Harish Chandra Solutions to
Macdonald's difference equations}
\label{sec4}
\begin{theorem}(Harish Chandra Solutions as matrix elements)
Let k be real and $0 < k < 1$.
1. We set  $r=  \frac{1}{1-k}$,$q=x^{2r}$, and  recall that $t=q^k$.

 For $\mu ={\lambda}\frac{r}{r-1}$
 the following matrix elements \footnotemark{}
\begin{eqnarray}
\phi(w \lambda+\rho, z)=  \prod z_i^{  k (\omega_1,\delta)
-\frac{r-1}{2 r}(\omega_1, \omega_1)}
\prod z_i^{-(i-1)\frac{r-1}{r}(\omega_1, \omega_1)}
\prod_{i < j} g_1\left(\frac{z_i}{z_j}\right)^{-1} \nonumber \\ \times
<0,\mu^*| \Phi_{h_{w(n)}}(z_n)\ldots  \Phi_{h_{w(2)}}(z_2)   
\Phi_{h_{w(1)}}(z_1)|0, \mu>\nonumber \\= (-1)^{\frac{(n-1)n}{2}}
\prod_{\alpha \in \Sigma_{+}} q^{\frac{(w \lambda, \alpha)((w \lambda, \alpha)+k)}{2}}
 \frac{\Gamma_q(1-k)}{\Gamma_q( (w \lambda, \alpha)+1) 
\Gamma_q(-( w\lambda, \alpha) +1-k)} z^{w \lambda +\rho}(1 +\ldots)
\label{coeff}
\end{eqnarray}
are the Harish Chandra  solutions to Macdonald's system of  difference equations (\ref{syst1}).
These solutions are correctly defined and converge to analytic
functions at least in the asymptotic zone $0 <|z_1| <|z_2| <\ldots <  |z_n|$. 

The function  $g_1(z)$ is defined by (\ref{g1}). 

\footnotetext{Recall the modification for the basic operators and formula
(\ref{feq1})}
2. Set $r=\frac{1}{k}$. For $\mu =\lambda \frac{r}{r-1}$ the following
matrix elements  (suitably modified)

\begin{eqnarray}
\prod z_i^{k(\omega_1,\delta)}
\prod_{i<j} z_j^{1-2k} 
\frac{(q^{k}\frac{z_i}{z_j};q)_{\infty}}{(q^{1-k} \frac{z_i}{z_j};q)_{\infty}}
\prod
z_i^{-\frac{r-1}{2r}(\omega_1, \omega_1)}
\prod z_i^{-(i-1)\frac{r-1}{r}(\omega_1, \omega_1)}
\prod_{i < j} g_1\left(\frac{z_i}{z_j}\right)^{-1} \nonumber \\ \times
<(0,\mu)^*| \Phi_{h_{w(n)}}(z_n)\ldots  \Phi_{h_{w(2)}}(z_2)   
\Phi_{h_{w(1)}}(z_1)|0, \mu>\nonumber \\= (-1)^{\frac{(n-1)n}{2}}
\prod_{\alpha \in \Sigma_{+}} q^{\frac{(w \lambda, \alpha)((w \lambda, \alpha)+1-k)}{2}}
 \frac{\Gamma_q(k)}{\Gamma_q( (w \lambda, \alpha)+1) 
\Gamma_q(-( w\lambda, \alpha) +k)} z^{w \lambda +\rho}(1 +\ldots)
\label{coeff1}
\end{eqnarray}
satisfy Macdonald's difference equations and produce solutions with
leading asymptotic $z^{w \lambda +\rho}$.


\label{th1}
\end{theorem}
\begin{remark}
The integral representations in theorem are the generalizations to the
Macdonald's case of integral representations for Heckman-Opdam
hypergeometric functions . In   the case of $A_1$ for the usual Gaussian
hypergeometric function $ z_1^{\lambda_1+\frac{k}{2}} 
z_2^{\lambda_2-\frac{k}{2}}F(k, \lambda_1-\lambda_2+k,
\lambda_1-\lambda_2+1, \frac{z_1}{z_2})$ 
one has the following integral representations

\begin{equation}
(z_1 z_2)^{\lambda_1+\frac{k}{2}} \int t^{\lambda_2-\lambda_1+k}
(z_1-t)^{-k} (z_2-t)^{-k} \frac{dt}{t}
\label{}
\end{equation}

and 
\begin{equation}
(z_1 z_2)^{\lambda_1+\frac{k}{2}} \int
t^{\lambda_2-\lambda_1-k}(z_1-t)^{k-1}(z_2 -t)^{k-1} dt 
\label{}
\end{equation}
correspondingly.
\end{remark}

\begin{remark}

Note  if $w \lambda+\rho$ is nonpositive integral weight,
then the Harish Chandra solution $\phi( w \lambda+\rho,z)$  times $\prod
z_i^{-(\lambda +\rho, \omega_1)}$: \quad\quad
$\prod z_i^{-(\lambda +\rho, \omega_1)} \phi( w \lambda+\rho,z)$
becomes Macdonald polynomial.
Compare also with a particular case of \cite{AKOS1} formulas $(35)$, $(37)$
with  $r^a=n-1-a$, $a=1, \ldots, n-1$, and allow parameters
$s^a$ to be nonintegral.  Parameters $r^a, s^a$ of  ref.
\cite{AKOS1} are integral and   related to the
Young diagram. Allowing parameters $s_a$ to be nonintegral , integral
representation of \cite{AKOS1} would rather give analogue of zonal
spherical function, i.e. solution to the system of difference
equations (\ref{syst1}), which admits continuation to single-valued
analytic function at $z_1=z_2=\ldots=z_n$. The type and the number of
variables $y_{ij}$ of integration that we use is rather remembrance of
the flag manifold and is minimal for generic parameter $\lambda$. 

 We indicate that parameter $k$ of this paper
corresponds to parameter $\beta$ of refs. \cite{AKOS1},\cite{AKOS2},\cite{FF}. 
\vskip 0.2 cm
We refer also to the results of ref. \cite{MN}.
\vskip 0.2 cm
Also, the situation is a deformation to the Macdonald's case of
integral representation for Heckman-Opdam hypergeometric functions
of theorem 3.2 in ref. \cite{A3}.
\end{remark}

\begin{theorem}(Analytic continuation) Consider
analytic continuation from the asymptotic zone :
$$0< |z_1| < \ldots <  |z_i| <  |z_{i+1}| <\ldots < |z_n|
$$
to the 
asymptotic zone 
$$
0 < |z_1|< \ldots <  |z_{i+1}| < |z_{i}| < \ldots < |z_n|
$$
such that
$$
arg\left(- \frac{z_i}{z_{i+1}} \right) < \pi.
$$
Let $\sigma_{i}$ permutes $i$th  and $i+1$th coordinate:
$$\sigma_i(z_1, z_2, \ldots, z_i, z_{i+1}, \ldots, z_n)=
(z_1, z_2, \ldots, z_{i+1},z_i,\ldots, z_n)
$$ 

Then analytic continuation of Harish Chandra solution to Macdonald's
difference equation $\phi(w \lambda+\rho, z)$ is given by: 
\begin{eqnarray}
\phi(w \lambda+\rho, z)\nonumber \\
=\frac{\Theta_q(q^k)}{\Theta_q(q^{\lambda_{w(i+1)}-\lambda_{w(i)}})}
\frac{\Theta_q(q^{\lambda_{w(i+1)}-\lambda_{w(i)}}
\frac{z_{i+1}}{z_{i}})}{\Theta_q(q^{k}\frac{z_{i+1}}{z_{i}})}
\left(\frac{z_i}{z_{i+1}}\right)^{\lambda_{w(i)}-\lambda_{w(i+1)}+k} 
{\phi(w \lambda+\rho,\sigma_i(z))}+
\nonumber \\
q^{k(\lambda_{w(i)}-\lambda_{w(i+1)})}
\frac{\Theta_q(q^{\lambda_{w(i)}-\lambda_{w(i+1)}+k})
}{\Theta_q(q^{\lambda_{w(i)}
-\lambda_{w(i+1)}})}
\frac{\Theta_q(\frac{z_{i+1}}{z_{i}})}{\Theta_q(q^{k}\frac{z_{i+1}}{z_{i}})}
\left(\frac{z_i}{z_{i+1}}\right)^{k}
     \phi(\sigma_i w\lambda+\rho,\sigma_i(z))
\label{f21}
\end{eqnarray}

\label{th2}
\end{theorem}
{\bf Proof:}
The theorem  \ref{th2} immediately follows from braiding properties of
 vertex operators $\Phi_{h_{w(i)}}(z_i)$ and
$\Phi_{h_{w(i+1)}}(z_{i+1})$ (\ref{f19}).
Note also that  if $w(i) < w(i+1)$ then since $w$ is a permutation
$$
< \mu +h_{w(1)} + \ldots +h_{w(i-1)}, \alpha_{w(i)} + \ldots
+\alpha_{w(i+1)-1}>=< \mu , \alpha_{w(i)} + \ldots+\alpha_{w(i+1)-1}>
$$ 

Thus the theorem follows. See also proposition \ref{prop4} below. 
$$\Box$$
\begin{corollary} Let $k$ be real and $0<k<1$.
If parameters $\lambda$ satisfies the following nondegeneracy conditions

 for any $l_1, l_2, \ldots, l_n \in \Bbb Z$ 
s.t. $\sum l_i^2 \ne 0$ one has 
\begin{equation}
\sum_{i=1}^n q^{\lambda_i + l_i} - \sum_{i=1}^n q^{\lambda_i} \ne 0
\label{res2} 
\end{equation}

Then there exists and  unique analytic solution to the Macdonald's system of difference
equations (\ref{Mac}) 
$$
D_z^m f( w\lambda +\rho, z)= c_{\lambda+\rho}^m f( w \lambda +\rho, z)
$$
$$m=1, \ldots,n$$
in the asymptotic zone 
$$0< |z_1 | < |z_2 | < \ldots < | z_n| $$
which has  the following  series expansion 
\begin{equation}
f(w \lambda+\rho, z)=z^{w \lambda+\rho}\sum_{p_1 , p_2 
,\ldots , p_{n-1} \in \Bbb Z_{+}} a(p_1, p_2, \ldots,p_{n-1})
\left(\frac{z_1}{z_2}\right)^{p_1} \ldots \left(\frac{z_{n-1}}{z_n}\right)^{p_{n-1}}
\label{ser1}
\end{equation}
\label{cor1}
with leading asymptotic coefficient equal to one:
 $a(0,0,\ldots,0)=1$.
The series expansion (\ref{ser1})  converges for sure  at least for 
$$0< |z_1 | < |z_2 | < \ldots < | z_n|. $$
\end{corollary}
{\bf Proof:}
In fact Harish Chandra  solution $\phi(w \lambda+\rho,z)$ is provided by   multiple integral of some
meromorphic function with compact cycle of integration. As long as we
can correctly draw the contours of integration it is correctly
defined, finite and analytic.
  Condition (\ref{res2}) guarantees uniqueness and along with
commutativity and special form of difference operators (\ref{f1}) that the
solution to the first order difference equation with leading exponent
$z^{w \lambda +\rho}$ is also a solution 
to the whole system of Macdonald's difference equations (\ref{syst1}).
We have two integral representations which give different restrictions
on parameter $\lambda$. Namely, we integrate analytic function over
compact cycle, as long as the contours for integrations enclose the
necessary poles and  the leading asymptotic coefficient is not equal
to zero, so that we can divide by it, we get convergence.
 Now we use the Hartogs theorem that analytic
function in $(\mathbb C^2 \setminus 0)$ can be extended to analytic function
on $\mathbb C^2$,
so only the condition $$\prod_{i<j}
\frac{1}{\Gamma_q(\lambda_{w(i)}-\lambda_{w(j)} +1)}\ne 0$$
comes from leading asymptotic coefficient , which is included in 
condition (\ref{res2}). 
$\Box$

\smallskip
{\bf Proof of theorem \ref{th1} :}

We proceed as follows.  Using Wick's rule rewrite the matrix elements as
multiple integrals with usual integration. Then we assume first that
$k < -1$ and $$ 0<<|z_1|<<|z_2|<<\ldots <<|z_{n+1}|.$$
Using propositions \ref{prop1}, \ref{prop2}, \ref{prop3} below we
prove by induction that integrals do satisfy the first order
difference equation. We use that the action of shift operators in the
above hypotheses does not change the disposition of contours of
integration and poles. 
The considered matrix elements permit series expansion (\ref{ser1})
and thus the coeficients of this series expansion satisfy certain
recurrent relations. The coefficients are essentially the finite sum
of ratios q-gamma functions and thus recurrent relations are satified
if parameter $k$ is extended from $k<-1$ to $k<1$ and then the domain 
$ 0<<|z_1|<<|z_2|<<\ldots <<|z_{n+1}|$ 
is replaced by $ 0<|z_1|<|z_2|<\ldots <|z_{n+1}|$ restricting to $0<k<1$. 
 Part 2 is proved using proposition \ref{propp} below.

\begin{proposition}(Change of variables in the  usual  integral)
\begin{equation}
\oint_{\gamma} (T_{q^{-1}} f(y)) h(y) dy=
q \oint_{q \gamma}  f(y) h(qy) d y=
q \oint_{q \gamma} f(y) (T_q h(y)) dy  
\label{f22}
\end{equation}
\label{prop1}
\end{proposition}
$\Box$

\begin{proposition}
Let $z=(z_1, z_2, \ldots, z_{n+1})$ , $y=(y_1,y_2, \ldots, y_n)$.

Let also $\Pi(z,y)$ be 
\begin{equation}
\Pi(z,y)= \prod_{i,j} \frac{ (q^{\frac{1+k}{2}}
\frac{z_i}{y_j};q)_{\infty}}
{(q^{\frac{1-k}{2}} \frac{z_i}{y_j};q)_{\infty}}
\label{f23}
\end{equation}

Then
\begin{equation}
D_z^1 (q,t) \Pi(z,y)= t \left(t^{n-1} D_y^1(q^{-1},t^{-1}) +1 \right)\Pi(z,y)
\label{f24}
\end{equation}
where $D_z^1$ $(D_y^1)$ is a first order difference operator defined
by (\ref{f1}).

\label{prop2}
\end{proposition}

Proposition immediately follows from Macdonald \cite{Ma1}
VI.2.13 on  page 311, and page 321.

         $\Box$

\begin{proposition}
\begin{eqnarray}
\left[\sum_{i} T_{q,y_i} \prod_{j \ne i} \frac{
\frac{1}{t}y_i-y_j}{y_i-y_j}\right]
\prod_{l < s} \frac{(\frac{y_l}{y_s};q)_{\infty} 
(\frac{y_s}{y_l};q)_{\infty}}
{(q^k \frac{ y_l}{y_s};q)_{\infty}(q^k \frac{y_s}{y_l};q)_{\infty}}
\nonumber \\=
\prod_{l <s}  \frac{(\frac{y_l}{y_s};q)_{\infty}
(\frac{y_s}{y_l};q)_{\infty}}
{(q^k \frac{y_l}{y_s};q)_{\infty} (q^k \frac{y_s}{y_l};q)_{\infty}}
 t^{-n+1} \sum_i \prod_{j \ne i} \frac{t y_i-y_j}{y_i-y_j} T_{q,y_i}
\label{f25}
\end{eqnarray}
\label{prop3}
\end{proposition}
 $\Box$

\begin{proposition}
\begin{eqnarray}
 \left[\sum_i \prod_{j \ne i} \frac{t z_i -z_j}{z_i-z_j} T_{q,z_i} \right]
 \prod_{i<j} z_j^{1-2k} \frac{(q^k \frac{z_i}{z_j};q)_{\infty}}
{(q^{1-k} \frac{z_i}{z_j};q)_{\infty}}\nonumber \\=
 \prod_{i<j} z_j^{1-2k} \frac{(q^k \frac{z_i}{z_j};q)_{\infty}}
{(q^{1-k} \frac{z_i}{z_j};q)_{\infty}}
\sum_i \prod_{j \ne i} \frac{\frac{q}{t} z_i -z_j}{z_i-z_j} T_{q,z_i}
\end{eqnarray}
\label{propp}
\end{proposition}
$\Box$
\begin{remark}
The above proposition gives rise to the transformation formulas
analogous to the following one for q-hypergeometric series:
\begin{eqnarray}
F_q(k, \lambda_1-\lambda_2+k, \lambda_1-\lambda_2+1, q^{1-k}
z)\nonumber \\=
\frac{(q^k z;q)_{\infty}}
{(q^{1-k} z;q)_{\infty}}
F_q(\lambda_1-\lambda_2+1-k, 1-k, \lambda_1-\lambda_2+1, q^k z)
\end{eqnarray}
\end{remark}

Here for the convenience of the reader,
 we recall the contraction formulas , cf. \cite{AJMP} , appendix C.

\begin{eqnarray}
:U_{\omega_1}(z_2): :U_{\omega_1}(z_1):= z_2^{\frac{r-1}{r}\frac{n-1}{n}}
g_1(\frac{z_1}{z_2}): U_{\omega_1}(z_2) U_{\omega_1}(z_1): \\
:U_{-\alpha_j}(z_1): : U_{\omega_j}(z_2): = z_1^{-\frac{r-1}{r}} s(\frac{z_2}{z_1})
:U_{-\alpha_j}(z_1) U_{\omega_j}(z_2): \\
:U_{\omega_j}(z_2): : U_{-\alpha_j}(z_1):
=z_2^{-\frac{r-1}{r}} 
s(\frac{z_1}{z_2}):U_{\omega_j}(z_2) U_{-\alpha_j}(z_1):\\
:U_{-\alpha_j}(z_1): :U_{-\alpha_{j+1}}(z_2):
 = z_1^{-\frac{r-1}{r}}
s(\frac{z_2}{z_1}) :U_{-\alpha_j}(z_1) U_{-\alpha_{j+1}}(z_2):\\
:U_{-\alpha_{j+1}}(z_2): :U_{-\alpha_j}(z_1): =z_2^{-\frac{r-1}{r}}
s(\frac{z_1}{z_2}) : U_{-\alpha_j}(z_1) U_{-\alpha_{j+1}}(z_2):\\
:U_{-\alpha_j}(z_1): : U_{-\alpha_j}(z_2): =z_1^{2 \frac{r-1}{r}}
t(\frac{z_2}{z_1}): U_{-\alpha_j}(z_1) U_{-\alpha_j}(z_2):
\label{f26} 
\end{eqnarray}

\begin{eqnarray}
s(z)= \frac{(x^{2r-1}z; x^{2r})_{\infty}}{(xz;x^{2r})_{\infty} }=
\frac{(q^{\frac{1+k}{2}}z;q)_{\infty}}
{(q^{\frac{1-k}{2}} z; q)_{\infty}} \\
t(z) = (1-z) \frac{(x^2 z; x^{2r})_{\infty}}{(x^{2r-2}z;
x^{2r})_{\infty}}=
(1-z) \frac{(q^{1-k} z; q)_{\infty}}
{(q^k z; q)_{\infty}}
\label{f27} 
\end{eqnarray}


\begin{eqnarray}
:U_{\omega_1}(z): : U_{-\alpha_1}(y):=
\frac{(q^{\frac{1+k}{2}} \frac{y}{z};q)_{\infty}}
{(q^{\frac{1-k}{2} } \frac{y}{z};q)_{\infty}} z^{-k}
:U_{\omega_1}(z)U_{-\alpha_1}(y): \nonumber \\
\frac{(q^{\frac{1+k}{2}} \frac{z}{y};q)_{\infty}}
{(q^{\frac{1-k}{2}} \frac{z}{y};q)_{\infty} }
z^{-k} \frac{\Theta_q(q^{\frac{1+k}{2}} \frac{y}{z})}
{\Theta_q(q^{\frac{1-k}{2}} \frac{y}{z})}:U_{\omega_1}(z)U_{-\alpha_1}(y):
\label{f28}
\end{eqnarray}

Under the action of shifts operators 
$\frac{\Theta_q(q^{\frac{1+k}{2}} \frac{y}{z})}
{\Theta_q(q^{\frac{1-k}{2}} \frac{y}{z})} z^{-k}$
behaves as $y^{-k}$.

\begin{eqnarray}
:U_{-\alpha_1}(y): : U_{\omega_1}(z):=
y^{-k} 
\frac{(q^{\frac{1+k}{2}} \frac{z}{y};q)_{\infty}}
{(q^{\frac{1-k}{2}} \frac{z}{y};q)_{\infty}}:U_{-\alpha_1}(y)
U_{\omega_1}(z):     
\label{f29}
\end{eqnarray}

\begin{eqnarray}
:U_{-\alpha_1}(y_1):  :U_{-\alpha_1}(y_2):=
y_1^{2k} \frac{(1-\frac{y_2}{y_1})(q^{1-k}
\frac{y_2}{y_1};q)_{\infty}}
{(q^k \frac{y_2}{y_1};q)_{\infty}}
:U_{-\alpha_1}(y_1) U_{-\alpha_1}(y_2):\nonumber \\
=
\frac{(\frac{y_1}{y_2};q)_{\infty} (\frac{y_2}{y_1};q)_{\infty}}
{(q^k \frac{y_1}{y_2};q)_{\infty}(q^k \frac{y_2}{y_1};q)_{\infty}}
y_1^{2k} \frac{\Theta_q(q^k \frac{y_1}{y_2})}
{\Theta_q(\frac{y_1}{y_2})} :U_{-\alpha_1}(y_1)U_{-\alpha_1}(y_2):
\label{f30}
\end{eqnarray}

Note, that term $y_1^{2k} \frac{\Theta_q(q^k \frac{y_1}{y_2})}
{\Theta_q(\frac{y_1}{y_2})}$ behaves under the action of shift operators 
as $(y_1 y_2)^k$.

We used:
$$
T_{q,z} \Theta_q(\alpha z)=
-\frac{1}{\alpha z}\Theta_q(\alpha z)
$$

$$T_{q,z} \frac{\Theta_q(\alpha z)}{\Theta_q(\beta z)}
=\frac{\beta}{\alpha}
\frac{\Theta_q(\alpha z)}{\Theta_q(\beta z)}
$$
and

$$
T_{q^{-1},z}
\Theta_q(\alpha z; q)_{\infty}=
-\frac{\alpha z}{
q} \Theta_q(\alpha z)
$$

$$
T_{q^{-1},z} \frac{\Theta_q(\alpha z)}{\Theta_q(\beta z)}=
\frac{\alpha}{\beta} \frac{\Theta_q(\alpha z)}{\Theta_q(\beta z)}
$$



Note
\begin{equation}
T_{q, x^{2 v}} [v]= [v+r]=-[v]
\label{f31}
\end{equation}
Thus the theorem  follows by induction. The base of induction is
checked directly.

\begin{eqnarray}(z_1 z_2)^{\frac{k}{2}} (z_1 z_2)^{-\frac{k(\omega_1, \omega_1)}{2}}
g_1(\frac{z_1}{ z_2})^{-1}
z_2^{-k(\omega_1,\omega_1)}
< (0, \mu +h_1 +h_2)^*| 
\Phi_{h_1}(z_2)
\Phi_{h_2}(z_1)
|0, \mu>\nonumber \\
=
(z_1 z_2)^{\lambda_1+\frac{k}{2}}
\int \frac{d y}{2 \pi i y}
y^{\lambda_2-\lambda_1+k}
z_1^{-k}
\frac{(q^{\frac{1+k}{2}}\frac{y}{z_1};q)_{\infty}}
{(q^{\frac{1-k}{2}}\frac{y}{z_1};q)_{\infty}}
z_2^{-k} \frac{(q^{\frac{1+k}{2}} \frac{y}{z_2};q)_{\infty}}
{(q^{\frac{1-k}{2}}\frac{y}{z_2};q)_{\infty}}\nonumber  \\
\times \frac{ \Theta_q(q^{\lambda_1-\lambda_2+\frac{1-k}{2}} \frac{y}{z_1})}
{\Theta_q(q^{\frac{k-1}{2}} \frac{y}{z_1})} 
( \frac{y}{z_1})^{\lambda_1-\lambda_2+1-k}
q^{\frac{1}{2}(\lambda_2-\lambda_1)(\lambda_2-\lambda_1+k) +\frac{k-1}{2}}
\nonumber \\
=
(z_1 z_2)^{\lambda_1+\frac{k}{2}} 
\int \frac{d y}{2 \pi i y}
y^{\lambda_2-\lambda_1+k}  
z_1^{-k}
\frac{(q^{\frac{1+k}{2}}\frac{y}{z_1};q)_{\infty}}
{(q^{\frac{1-k}{2}}\frac{y}{z_1};q)_{\infty}}
z_2^{-k} \frac{(q^{\frac{1+k}{2}} \frac{y}{z_2};q)_{\infty}}
{(q^{\frac{1-k}{2}}\frac{y}{z_2};q)_{\infty}}\nonumber  \\
\times \frac{ \Theta_q(q^{\lambda_1-\lambda_2+\frac{1-k}{2}} \frac{y}{z_1})}
{\Theta_q(q^{\frac{1+k}{2}} \frac{y}{z_1})} 
( \frac{z_1}{y})^{\lambda_2-\lambda_1+k}
(-1)q^{\frac{1}{2}(\lambda_2-\lambda_1)(\lambda_2-\lambda_1+k)}\nonumber \\=
(-1) q^{\frac{1}{2}(\lambda_2-\lambda_1)(\lambda_2-\lambda_1+k)}
\frac{\Gamma_q(1-k)}{\Gamma_q(\lambda_1-\lambda_2+1-k)
\Gamma_q(\lambda_2-\lambda_1+1)}\nonumber \\
\times z_1^{\lambda_2+\frac{k}{2}} z_2^{\lambda_1-\frac{k}{2}} F_q( k,
\lambda_2-\lambda_1+k, \lambda_2-\lambda_1+1, q^{1-k} \frac{z_1}{z_2})
\end{eqnarray}

\begin{eqnarray}
(z_1 z_2)^{\frac{k}{2}} 
(z_1 z_2)^{-\frac{k(\omega_1,\omega_1)}{2}}
z_2^{-k(\omega_1,\omega_1)}
g_1(\frac{z_1}{z_2})^{-1}
< (0, \mu +h_1 +h_2)^*| 
\Phi_{h_2}(z_2)
\Phi_{h_1}(z_1)
|0, \mu>\nonumber\\ 
(z_1 z_2)^{\lambda_1+\frac{k}{2}} 
\int \frac{d y}{2 \pi i y}
y^{\lambda_2-\lambda_1+k}
y^{-k}
\frac{(q^{\frac{1+k}{2}}\frac{z_1}{y};q)_{\infty}}
{(q^{\frac{1-k}{2}}\frac{z_1}{y};q)_{\infty}}
z_2^{-k} \frac{(q^{\frac{1+k}{2}} \frac{y}{z_2};q)_{\infty}}
{(q^{\frac{1-k}{2}}\frac{y}{z_2};q)_{\infty}}\nonumber \\
\times \frac{ \Theta_q(q^{\lambda_1-\lambda_2+\frac{1+k}{2}} \frac{y}{z_2})}
{\Theta_q(q^{\frac{k-1}{2}} \frac{y}{z_2})}
( \frac{y}{z_2})^{1+\lambda_1-\lambda_2}
q^{\frac{(\lambda_2-\lambda_1)(\lambda_2-\lambda_1-k)}{2} +\frac{k-1}{2}}
\nonumber \\=
(z_1 z_2)^{\lambda_1+\frac{k}{2}} 
\int \frac{d y}{2 \pi i y}
y^{\lambda_2-\lambda_1+k}  
y^{-k}
\frac{(q^{\frac{1+k}{2}}\frac{z_1}{y};q)_{\infty}}
{(q^{\frac{1-k}{2}}\frac{z_1}{y};q)_{\infty}}
z_2^{-k} \frac{(q^{\frac{1+k}{2}} \frac{y}{z_2};q)_{\infty}}
{(q^{\frac{1-k}{2}}\frac{y}{z_2};q)_{\infty}}
\nonumber \\
\times \frac{ \Theta_q(q^{\lambda_1-\lambda_2+\frac{1+k}{2}} \frac{y}{z_2})}
{\Theta_q(q^{\frac{1+k}{2}} \frac{y}{z_2})} 
( \frac{z_2}{y})^{\lambda_2-\lambda_1}
(-1)q^{\frac{(\lambda_2-\lambda_1)(\lambda_2-\lambda_1-k)}{2}}\nonumber \\=
(-1) q^{\frac{(\lambda_2-\lambda_1)(\lambda_2-\lambda_1-k)}{2}}
\frac{\Gamma_q(1-k)}{\Gamma_q(\lambda_2-\lambda_1+1-k)
\Gamma_q(\lambda_1-\lambda_2+1)}\nonumber  \\
\times z_1^{\lambda_1+\frac{k}{2}} z_2^{\lambda_2-\frac{k}{2}} F_q( k,
\lambda_1-\lambda_2+k, \lambda_1-\lambda_2+1, q^{1-k} \frac{z_1}{z_2})
\end{eqnarray}

Finally , we note that Harish Chandra solution $\phi(w \lambda +\rho, z)$
being an eigenfunction of the first order difference equation $D^1_z $
is also a solution to the whole system of difference equations
(\ref{Mac}) if the following condition holds:

for any $l_1, l_2, \ldots, l_n \in \Bbb Z$ 
s.t. $\sum l_i^2 \ne 0$ one has 
\begin{equation}
\sum_{i=1}^n q^{\lambda_i + l_i} -\sum_{i=1}^n   q^{\lambda_i} \ne 0
\label{res1} 
\end{equation}

The leading asymptotic coefficient is equal:
\begin{eqnarray}
<(0, \mu +h_{w(1)})^*|\Phi_{h_{w(1)}}(1)|0,\mu>
<(0, \mu +h_{w(1)}+h_{w(2)})^*|\Phi_{h_{w(2)}}(1)|0,\mu+h_{w(1)}> \times
\ldots \nonumber \\ \times 
<(0, \mu +h_{w(1)}+h_{w(2)}+\ldots +h_{w(n)})^*|\Phi_{h_{w(n)}}(1)|0,\mu+h_{w(1)}+\ldots +
h_{w(n-1)}>
\label{}
\end{eqnarray}
and can be easily calculated using proposition \ref{prop} below.
The theorem is proved.
$\Box$

\begin{proposition}
\label{prop4} Recall that $\mu =\frac{\lambda}{k}$, $r=\frac{1}{1-k}$.
 For $i \ne j$ one has
\begin{eqnarray}
(z_1 z_2)^{\frac{k}{2}} (z_1 z_2)^{-\frac{k (\omega_1, \omega_1)}{2}}
z_2^{-k(\omega_1, \omega_1)}
g_1(\frac{z_1}{z_2})^{-1}     
<(0,\mu +h_i+h_j)^*| \Phi_{h_j}(z_2) \Phi_{h_i}(z_1) |0, \mu>
\label{s1} \\=
(-1)^{i+j-2} \prod_{s <i\; s \ne j} \frac{\Gamma_q(1-k)
q^{\frac{(\lambda_i-\lambda_s)(\lambda_i-\lambda_s+k)}{2}}}
{\Gamma_q(\lambda_i-\lambda_s+1) \Gamma_q(\lambda_s-\lambda_i+1-k)}
\nonumber \\
\prod_{s <j\; s \ne i} \frac{ \Gamma_q(1-k)
q^{\frac{(\lambda_j-\lambda_s)(\lambda_j-\lambda_s+k)}{2}}}
{\Gamma_q(\lambda_j-\lambda_s+1)
\Gamma_q(\lambda_s-\lambda_j+1-k)
}
\nonumber \\
\times \frac{\Gamma_q(1-k)
q^{\frac{(\lambda_i-\lambda_j)(\lambda_i-\lambda_j+k)}{2}}}
{\Gamma_q(\lambda_i-\lambda_j+1)
\Gamma_q(\lambda_j-\lambda_i+1-k)}\nonumber \\
\times z_1^{\lambda_i+\frac{k}{2}}
z_2^{\lambda_j-\frac{k}{2}}
 F_q(k,  \lambda_i-\lambda_j+k,  \lambda_i -\lambda_j+1, q^{1-k}\frac{z_1}{z_2})   
\end{eqnarray}

\end{proposition} 

{\bf Proof:}
In fact, matrix element (\ref{s1})
 satisfies Macdonald's difference equation (\ref{f4})
 (the proof is the
same as in theorem \ref{th1})  and has leading
exponent $z_1^{\lambda_i+\frac{k}{2}}z_2^{\lambda_j-\frac{k}{2}}$.
Thus it is equal to the q-hypergeometric function
$$z_1^{\lambda_i+\frac{k}{2}}z_2^{\lambda_j-\frac{k}{2}}
F_q(k, \lambda_i-\lambda_j+k, \lambda_i-\lambda_j+1, q^{1-k} \frac{z_1}{z_2})
$$
up to the leading asymptotic coefficient. The leading asymptotic
coefficient is calculated using proposition \ref{prop}.

The proposition is in complete agreement with 
 braiding properties 
(commutation relation) of  vertex operators $\Phi_{h_j}(z_2)$ ,$\Phi_{h_i}(z_1)$
 (\ref{f19}).

$$\Box$$

\begin{proposition}
Recall that $\mu =\frac{\lambda}{k}$, $r=\frac{1}{1-k}$. Let $i= 2,\ldots,n$.
\begin{eqnarray}
< (0, \mu + h_i)^*|\Phi_{h_i}(1) |0, \mu>=
(-1)^{i-1} \prod_{s=1}^{i-1} \frac{\Gamma_q(1-k)
q^{\frac{(\lambda_i-\lambda_s)(\lambda_i -\lambda_s+k)}{2}}}
{\Gamma_q(\lambda_i-\lambda_s+1)\Gamma_q(\lambda_s-\lambda_i +1-k)}
\label{}
\end{eqnarray}
\label{prop}
\end{proposition}
$\Box$

\begin{proposition} For any $ n \in \mathbb N$
\label{prop6}
\begin{eqnarray}
\oint y^n \frac{\Theta_q(q^{\lambda_2-\lambda_1+\frac{k+1}{2}}
\frac{1}{y})}
{\Theta_q(q^{\frac{1+k}{2}}\frac{1}{y})}
\frac{(q^{\frac{1+k}{2}} \frac{1}{y};q)_{\infty}}
{(q^{\frac{1-k}{2}} \frac{1}{y};q)_{\infty}}
\frac{d y}{2 \pi i y}\nonumber \\
=
\frac{\Gamma_q(1-k)}{\Gamma_q(\lambda_2-\lambda_1+1)
\Gamma_q(\lambda_1-\lambda_2+1-k)}
\frac{\Gamma_q(\lambda_2-\lambda_1+k+n)}
{\Gamma_q(\lambda_2-\lambda_1+k)}
\frac{\Gamma_q(\lambda_2-\lambda_1+1)}
{\Gamma_q(\lambda_2-\lambda_1+n+1)}
\label{}
\end{eqnarray}
Here the integration contour is a simple closed curve around the
origin, which encloses the poles
$$y=q^{\frac{1-k}{2}+m}$$
$m=0,1,\ldots$.

\end{proposition}
{\bf Proof:}
Let's assume that $n \ge 0$. The case $n<0$ is considered similarly.
\begin{eqnarray}
\int y^n \frac{\Theta_q(q^{\lambda_2-\lambda_1+\frac{k+1}{2}}
\frac{1}{y})}
{\Theta_q(q^{\frac{1+k}{2}}\frac{1}{y})}
\frac{(q^{\frac{1+k}{2}} \frac{1}{y};q)_{\infty}}
{(q^{\frac{1-k}{2}} \frac{1}{y};q)_{\infty}}
\frac{d y}{2 \pi i y} 
=\int y^n \frac{(q^{\lambda_1-\lambda_2+\frac{1-k}{2}} y;q)_{\infty}}
{(q^{\frac{1-k}{2} }y;q)_{\infty}}
\frac{(q^{\lambda_2-\lambda_1+\frac{k+1}{2}} \frac{1}{y};q)_{\infty}}
{(q^{\frac{1-k}{2}} \frac{1}{y};q)_{\infty}} \frac{ d y}{ 2 \pi i
y} \nonumber \\ 
 =\sum_{m=0}^{\infty}
\frac{(1-q^{\lambda_1-\lambda_2}) \ldots
(1-q^{\lambda_1-\lambda_2+m-1})}
{(1-q) \ldots (1-q^m)} q^{\frac{1-k}{2}m} \int y^m y^n
\frac{(q^{\lambda_2-\lambda_1+\frac{k+1}{2} }\frac{1}{y};q)_{\infty}}
{(q^{\frac{1-k}{2}} \frac{1}{y};q)_{\infty}} \frac{d y}{ 2 \pi i
y}=\nonumber \\
\sum_{m=0}^{\infty}
\frac{(1-q^{\lambda_1-\lambda_2}) \ldots
(1-q^{\lambda_1-\lambda_2+m-1})}
{(1-q) \ldots (1-q^m)} q^{\frac{1-k}{2}m} 
\frac{(1-q^{\lambda_2-\lambda_1+k}) \ldots
(1-q^{\lambda_2-\lambda_1+k+m+n-1})}
{(1-q) \ldots (1-q^{m+n})} q^{\frac{1-k}{2} (m+n)} \nonumber \\=
q^{\frac{1-k}{2}n} \frac{ (1-q^{\lambda_2-\lambda_1+k})\ldots
(1-q^{\lambda_2-\lambda_1+k+n-1})}
{(1-q)\ldots (1-q^n)} F_q( \lambda_2-\lambda_1+k+n,
\lambda_1-\lambda_2, 1-n, q^{1-k}) \nonumber \\
q^{\frac{1-k}{2}n}
\frac{\Gamma_q(\lambda_2-\lambda_1+k+n)}{\Gamma_q(\lambda_2-\lambda_1+k)}
 \frac{\Gamma_q(1-k)}{\Gamma_q(\lambda_1-\lambda_2+1-k)
\Gamma_q(\lambda_2-\lambda_1+n+1)}\nonumber \\
\label{v}
\end{eqnarray} 

In the above claculation we used that 
$$
F_q(a,b,c,q^{c-a-b})=\frac{\Gamma_q(c) \Gamma_q(c-a-b)}{\Gamma_q(c-a)
\Gamma_q(c-b)}
$$
and 
\begin{eqnarray} 
\oint y^n \frac{(\alpha \frac{1}{y};q)_{\infty}}
{(\beta \frac{1}{y};q)_{\infty} } \frac{ dy}{2 \pi i y}
= 
\oint y^{-n} \frac{(\alpha y;q)_{\infty}}{(\beta y;q)_{\infty}} \frac{d
y}{ 2 \pi i y} \nonumber \\
=
\frac{(1- \frac{\alpha}{\beta}) \ldots (1-q^{n-1}
\frac{\alpha}{\beta})}
{(1-q) \ldots (1-q^n)} \beta^n
\label{}
\end{eqnarray}
Here we assumed that the contour of integration is a simple closed
curve around the origin which encloses the poles $ y= \beta q^m$,
$m=0, 1, \ldots$.

The  calculation below shows that 
\begin{equation}
\lim_{N \to \infty} 
\oint_{ q^N \gamma} y^n \frac{\Theta_q(q^{\lambda_2-\lambda_1+\frac{k+1}{2}}
\frac{1}{y})}
{\Theta_q(q^{\frac{1+k}{2}}\frac{1}{y})}
\frac{(q^{\frac{1+k}{2}} \frac{1}{y};q)_{\infty}}
{(q^{\frac{1-k}{2}} \frac{1}{y};q)_{\infty}}
\frac{d y}{2 \pi i y}=0
\label{conv}
\end{equation}

 If we count the residues that contribute to the integral are at
$y=q^{\frac{1-k}{2}+m}$, $m=0,1,\ldots$, we get 

\begin{eqnarray}
\oint y^n \frac{\Theta_q(q^{\lambda_2-\lambda_1+\frac{k+1}{2}}
\frac{1}{y})}
{\Theta_q(q^{\frac{1+k}{2}}\frac{1}{y})}
\frac{(q^{\frac{1+k}{2}} \frac{1}{y};q)_{\infty}}
{(q^{\frac{1-k}{2}} \frac{1}{y};q)_{\infty}}
\frac{d y}{2 \pi i y}\nonumber \\
=
q^{\frac{1-k}{2} n}\frac{\Theta_q(q^{\lambda_2-\lambda_1+k})}
{\Theta_q(q^k)} \sum_{m=0}^{\infty}
q^{(\lambda_2-\lambda_1)m} q^{n m} 
Res_{y= q^{\frac{1-k}{2}+m}}
\frac{(q^{\frac{1+k}{2}} \frac{1}{y};q)_{\infty}}
{(q^{\frac{1-k}{2}} \frac{1}{y};q)_{\infty}}\frac{1}{y}
\nonumber \\
=
\frac{\Theta_q(q^{\lambda_2-\lambda_1+k})}{\Theta_q(q^k)}
q^{\frac{1-k}{2}n} \frac{(q^k;q)_{\infty}}{(q;q)_{\infty}}
\sum_{m=0}^{\infty} q^{ k m} \frac{(1-q^{1-k})\ldots (1-q^{m-k})}
{(1-q)\ldots (1-q^m)}q^{(\lambda_2-\lambda_1+n)m} \nonumber \\
=
\frac{q^{\frac{1-k}{2}n} \Theta_q(q^{\lambda_2-\lambda_1+k})
(q^k;q)_{\infty}}
{\Theta_q(q^k) (q;q)_{\infty}}
\frac{(q^{1-k}q^{\lambda_2-\lambda_1+k+n};q)_{\infty}}
{(q^{\lambda_2-\lambda_1+k+n};q)_{\infty}}\nonumber \\
=
\frac{q^{\frac{1-k}{2}n} (q^{\lambda_2-\lambda_1+k};q)_{\infty}
(q^{\lambda_1-\lambda_2+1-k};q)_{\infty}}
{(q;q)_{\infty}(q^{1-k};q)_{\infty}}
\frac{(q^{\lambda_2-\lambda_1+n+1};q)_{\infty}}
{(q^{\lambda_2-\lambda_1+k+n};q)_{\infty}} \nonumber \\
=\frac{\Gamma_q(1-k)}{\Gamma_q(\lambda_2-\lambda_1+1) 
\Gamma_q(\lambda_1-\lambda_2+1-k)}
q^{\frac{1-k}{2}n} \frac{\Gamma_q(\lambda_2-\lambda_1+k+n)}
{\Gamma_q(\lambda_2-\lambda_1+k)}
\frac{\Gamma_q(\lambda_2-\lambda_1+1)}
{\Gamma_q(\lambda_2-\lambda_1+n+1)} 
\label{}
\end{eqnarray}
In the calculation we used binomial expansion formula:
\begin{equation}
1 + \frac{1-q^a}{1-q} x + \frac{(1-q^a)(1-q^{a+1})}{(1-q)(1-q^2)}x^2 +
\ldots = \frac{(q^a x;q)_{\infty}}{(x;q)_{\infty}}
\label{}
\end{equation}
Since we got the same value as (\ref{v}) 
, we proved the limit (\ref{conv}).
$$\Box$$
\begin{remark}
The above calculation should be compared with

\begin{eqnarray}
\int t^{\lambda_2-\lambda_1+k} (1-t)^{-k} \frac{d t}{ 2 \pi i t}=
\frac{\Gamma(\lambda_2-\lambda_1+k)
\Gamma(1-k)}{\Gamma(\lambda_2-\lambda_1+1)} \times
\frac{ e^{2 \pi i (\lambda_2-\lambda_1+k)}-1}{2 \pi i} \nonumber \\
=
e^{ \pi i (\lambda_2-\lambda_1+k)} \frac{\Gamma(1-k)}
{\Gamma(\lambda_2-\lambda_1+1) \Gamma(\lambda_1-\lambda_2+1-k)}
\label{}
\end{eqnarray}
where  the contour for integration for $t$ is a closed loop which goes
counterclockwise and  starts and ends at $1$. 
\end{remark}

\begin{proposition}
For $\mu = \frac{\lambda}{k}$ , $r=\frac{1}{1-k}$ the following equality holds
\begin{eqnarray}
\prod g_1(\frac{z_i}{z_j})^{-1} \prod z_i^{-(i-1)k(\omega_1,\omega_1) }
\prod z_i^{-\frac{k}{2}(\omega_1,\omega_1)} 
<(0,\mu+n h_i)^*| \Phi_{h_i}(z_n) \ldots \Phi_{h_i}(z_1)|0,
\mu>\nonumber \\=
C (z_1 \ldots z_n)^{\lambda_i-\lambda_1}
\label{}
\end{eqnarray}
The leading asymptotic coefficient $C$ is easily obtained using
proposition \ref{prop}.
\end{proposition}
$\Box$
 
{\bf Acknowledgements.} 
We are grateful to I. Gelfand for interesting discussions and for referring us to
refs. \cite{MN},\cite{Mi},\cite{Ch} , to S. Lukyanov for interesting discussions   
, to H. Awata for interesting remarks and
comments and references to \cite{AKOS1}, \cite{AKOS2}.


\begin{thebibliography}{99}

\bibitem{AJMP}
Y.Asai, M.Jimbo, T. Miwa, and  Y.Pugai
\newblock Bosonization of vertex operators for the $A^{(1)}_{n-1}$
face model 
\newblock J. Phys. A29 1996, 6595-6616 
\newblock  hep-th/9606095,  June 1996

\bibitem{Askey}
Askey R.
\newblock Beta integrals in Ramanujan's papers, his unpublished work
and further examples
\newblock in Ramanujan revisited, ed. Andrews G.E., Askey R.A., 
Berndt B.C., Ramathan K.G., Rankin R.A.
\newblock Proceedings of Centenary conference University of Illinois
at Urbana-Champaign, June 1-5, 1987
\newblock 561-590
\bibitem{AKOS2}
Awata H., Kubo H., Odake S., Shiraishi J.
\newblock A quantum deformation of the Virasoro algebra and
Macdonald symmetric functions
\newblock preprint q-alg 9507034


\bibitem{AKOS1}
Awata H., Kubo H., Odake S., and  Shiraishi J.
\newblock Quantum $W_n$ algebras and Macdonald polynomials 
\newblock {\em Comm. Math. Phys. 179} 401--416, 1996.

\bibitem{AKOS3}
Awata H., Kubo H., Odake S., and Shiraishi J.
\newblock Quantum deformation of the $W_n$ algebra
\newblock q-alg/9612001 

\bibitem{AKOS4}
Awata H., Odake S., and J. Shiraishi
\newblock Integral representations of the Macdonald Symmetric 
Functions
\newblock q-alg/9506006

\bibitem{}
Baxter R.J.
\newblock Exactly Solved Models in Statistical Mechanics
\newblock Academic Press, London, 1982


\bibitem{}
Bernstein J., Gelfand I., Gelfand S.
\newblock Differential operators on the base affine space and a
study of $\mathfrak g$-modules
\newblock in  Lie groups and their representations
\newblock{ Proc. Summer school in group reps. Bolyai Math. Soc., 
Budapest}  21--64, 1971. 

\bibitem{}
Bailey W.N.
\newblock Generalized hypergeometric series
\newblock {\em Cambridge University press}, 1935

\bibitem{Ch}
Cherednik I.
\newblock Difference Macdonald-Mehta conjecture
\newblock q-alg/9702022

\bibitem{Exton}
Exton H.
\newblock q-hypergeometric functions and applications
\newblock Ellis Horwood Limited, 1983


\bibitem{FL}
Fateev V.,  Lukyanov S.
\newblock The models of two-dimensional conformal quantum field theory 
with $Z_n$ symmetry.
\newblock{\em Int. J. mod.phys A } 3: 507--520, 1988

\bibitem{FL1}
 Fateev V.,  Lukyanov S.
Additional symmetries and exactly soluble models in two
dimensional field theory.
\newblock {\em Sov. sci. Rev.} 15:1--117, 1990

\bibitem{FZ}
 Fateev V.,  Zamolodchikov A.
\newblock Conformal quantum field theory models
in two dimensions having $\mathbb Z_3$ symmetry
\newblock{\em Nucl. Phys. B280} 644--660, 1987.

\bibitem{Fe}
Felder G.
\newblock BRST approach to minimal models.
\newblock {\em Nucl. Phys. B}, 317:215--236, 1989.

\bibitem{}
Felder G., Varchenko A.
\newblock Elliptic quantum groups and Ruijsenaars models
\newblock q-alg 9704005

\bibitem{FF}
Frenkel E., Feigin B.
\newblock Quantum W-algebras and elliptic algebras
\newblock {\em Comm. Math. Phys. A3}, 237--264, 1996.

\bibitem{}
Frenkel E., Reshetikhin N.
\newblock Quantum affine algebras and deformations of the Virasoro and
W-algebras
\newblock q-alg/9505025


\bibitem{}
Frenkel I., Reshetikhin N.
\newblock Quantum affine algebras and holonomic difference equations.
\newblock {\em Comm. math. phys. 145} 755--815, 1992.


\bibitem{Gel}
Gelfand I.
\newblock Spherical functions on symmetric Riemannian spaces.


\newblock {\em Dokl. Akad. Nauk USSR}, 70:5--8, 1950.
\bibitem{GGP}
 Gelfand I., Graev M., Piatetskii-Shapiro I.
\newblock Representation theory and automorphic functions
\newblock Generalized functions vol. 6


\bibitem{GN1}
Gelfand~I.M., Naimark~M.A. 
\newblock Unitary representations of classical groups.
\newblock {\em Tr.Mat.Inst. Steklova}, 36:1--288, 1950.


\bibitem{GinK}
Gindikin S.G., Karpelevich F.I.
\newblock Plancherel measure for Riemannian symmetric spaces 
of nonpositive curvature
\newblock{\em Dokl. Akad. Nauk USSR} 145(2):252--255,1962.




\bibitem{Hel}
Helgason S.
\newblock {\em Groups and geometric analysis}.
\newblock Academic Press, Inc., 1984.


\bibitem{Ha}
Harish Chandra
\newblock Spherical functions on a semi-simple Lie group I.
\newblock {\em Amer. J. of Math}, 80: 241--310, 1958.




\bibitem{HO}
Heckman G., Opdam E.
\newblock Root systems and hypergeometric functions I
\newblock {\em Comp. math. 64} 329--352, 1987.

\bibitem{}
Heckman G.
\newblock Root systems and hypergeometric functions II
\newblock {\em Comp. math. 64} 353--373, 1987.

\bibitem{}
Heckman G., Opdam E.
\newblock Yang's system of particles and Hecke algebras
\newblock{ \em Annals of Math 145} 139--173, 1997.

\bibitem{JLMP}
Jimbo M., Lashkevich M., Miwa T., Pugai Y.
\newblock Lukyanov's Screening operators for the deformed Virasoro
algebra
\newblock preprint hep-th 9607177, July 1996.

\bibitem{A1}
Kazarnovski-Krol A.
\newblock Value of generalized hypergeometric function at unity.
\newblock to appear in Arnold-Gelfand mathematical  seminars, Birkhauser,Boston, 1996



\bibitem{A2}
Kazarnovski-Krol A.
\newblock Harish Chandra decomposition for zonal spherical functions
of type $A_n$
\newblock to appear in Arnold-Gelfand mathematical seminars,
Birkhauser, Boston, 1996



\bibitem{A3}
Kazarnovski-Krol A.
\newblock A generalization of Selberg integral
\newblock q-alg 9507011

\bibitem{A4}
Kazarnovski-Krol A.
\newblock Cycles for the asymptotic solutions and the Weyl group
\newblock{\em in Gelfand mathematical seminars 1993-1995,
ed. by I.Gelfand, J.Lepowsky, M. Smirnov, Birkhauser}, 122--150, 1996

\bibitem{A5}
Kazarnovski-Krol A.
\newblock Cycle for integration for zonal spherical function of type
$A_n$
\newblock q-alg 9511008



\bibitem{L2}
Lukyanov S.
\newblock Quantization of Gelfand-Dikii bracket
\newblock{\em Funct.anal. appl.} 22:4, 1988

\bibitem{LP}
Lukyanov S., Pugay Y.
\newblock Multi-point Local height probabilities in the integrable
RSOS Model
\newblock hep-th 9602074
\newblock {\em Nucl. Phys B473}, 631--658, 1996

\bibitem{LP2}
Lukyanov S., Pugay Y.
\newblock Bosonization of ZF Algebras: Direction toward Deformed
Virasoro Algebra
\newblock hep-th 9412128


\bibitem{LL}
Lukyanov S.
\newblock A note on the deformed Virasoro algebra
\newblock hep-th 9509037

\bibitem{Ma1}
Macdonald I.G.
\newblock Symmetric functions and Hall polynomials
\newblock second edition Clarendon Press, Oxford University press, 1995 

\bibitem{Ma2}
Macdonald I.G.
\newblock Spherical functions on a group of p-adic type
\newblock{\em Publication of Ramanujan Institute, No. 2,  Madras}, 1971.

\bibitem{Ma3}
Macdonald I.G.
\newblock Commuting differential operators and zonal spherical
functions
\newblock{\em Springer lecture notes , 1271}, 189--200, 1987.



\bibitem{MN}
Mimachi K., Noumi M.
\newblock Notes on eigenfunctions for Macdonald's q-difference
operators
\newblock preprint 1996.

\bibitem{Mi}
Mimachi K.
\newblock A solution to quantum Knizhnik-Zamolodchikov equations and
its application to eigenvalue problems of the Macdonald type
\newblock to appear Duke math. J.


\bibitem{Sek}
Sekigushi J.
\newblock Zonal spherical functions on some symmetric spaces.
\newblock {\em Publ. RIMS, Kyoto University}, 12:455--459.


\bibitem{Za1}
Zamolodchikov A.B.
\newblock Infinite additional symmetries in two-dimensional
quantum field theory
\newblock{\em Teor.  math. phys} 65(3):1205--1213, 1986.

\end{thebibliography}
\end{document}